\input amstex
\documentstyle{amsppt}
\magnification=1200
\loadeusm
\loadeufb
\loadbold
\def\qed{\hfill$\ssize\square$}
\predefine\accute{\'}
\redefine\'{\kern.05em{}}
\def\<{\langle}
\def\>{\rangle}

\def\geqs{\geqslant}
\def\im{\operatorname{im}}
\def\Tr{\operatorname{Tr}}
\def\Skew{\operatorname{Skew}}
\def\Ric{\operatorname{Ric}}
\def\div{\operatorname{div}}
\def\Ad{\operatorname{Ad}}
\def\Nil{\text{\it Nil\/}}
\def\Sol{\text{\it Sol\/}}
\def\D{\Cal D}
\def\E{\Cal E}
\def\V{\Cal V}
\def\R{\eusm R}
\def\bD{{\fam\cmbsyfam D}}
\def\bE{{\fam\cmbsyfam E}}

\def\bfH{\eufb H}
\def\bfM{\eufb M}
\def\a{\alpha}
\def\b{\beta}
\def\s{\sigma}
\def\bPhi{\boldsymbol\Phi}
\def\bxi{\boldsymbol\xi}
\def\bJ{\boldsymbol J}
\def\boldeta{\boldsymbol\eta}
\def\fg{\frak g}
\def\fh{\frak h}
\def\fm{\frak m}
\def\fH{\frak H}
\def\fM{\frak M}
\def\nab#1{\nabla\kern-.2em\lower.8ex\hbox{$\ssize#1$}\'}
\def\nabv#1{\nabla^v\kern-.4em\lower.8ex\hbox{$\ssize#1$}\'}
\def\nabc#1{\nabla^c\kern-.4em\lower.8ex\hbox{$\ssize#1$}\'}
\def\nabsq#1#2{\nabla^2\kern-.3em\lower.8ex\hbox{$\ssize#1,#2$}\,} 
\def\midstar{\'\hbox{$*$}\'}
\def\uperp#1{#1\raise1.1ex\hbox{$\sssize\bot$}}
\def\dvs{d^{\'v}\kern-.15em\s}
\PSAMSFonts

\NoBlackBoxes
\nologo
\TagsOnRight
\rightheadtext{Harmonic Almost Contact Structures}
\pagewidth{16truecm}
\pageheight{24.3truecm}
\linespacing{1.3}
\topmatter
\title
Harmonic Almost Contact Structures
\endtitle
\author 
E\. Vergara-Diaz
\&\
C\. M\. Wood
\endauthor
\address 
Department of Mathematics, University of York, Heslington, York
Y010 5DD, UK.
\endaddress
\email 
evd103\@york.ac.uk,
cmw4\@york.ac.uk 
\endemail
\thanks
\endthanks
\abstract
An almost contact metric structure is parametrized by a section $\s$ of an
associated homogeneous fibre bundle, and conditions for $\s$ to be a
harmonic section, and a harmonic map, are studied.  These involve the
characteristic vector field $\xi$, and the almost complex structure in the
contact subbundle.  Several examples are given where the harmonic section
equations for $\s$ reduce to those for $\xi$, regarded as a section of the
unit tangent bundle.  These include trans-Sasakian structures, and certain nearly
cosymplectic structures.  On the other hand, we obtain examples where
$\xi$ is harmonic but $\s$ is not a harmonic section.  Many of our examples
arise by considering hypersurfaces of almost Hermitian manifolds, with the
induced almost contact structure, and comparing the harmonic section
equations for both structures.
\endabstract
\keywords
Harmonic section, harmonic map, harmonic unit vector field, almost contact
metric structure, contact metric structure, trans-Sasakian, nearly
cosymplectic, nearly K\"ahler structure
\endkeywords
\subjclass
53C10, 53C15, 53C43, 53C56, 53D10, 53D15, 58E20
\endsubjclass
\endtopmatter
\document
\head
1. Introduction
\endhead
An almost contact metric structure on an orientable manifold $M$ of odd
dimension $2n+1$ may be regarded as a reduction of the structure group of
$M$ to $U(n)$ \cite{15}.  However the holonomy of such a manifold lies in
$U(n)$ if and only if the almost contact structure is {\sl cosymplectic\/}
\cite{1}, which in particular precludes all contact metric structures. 
Thus, if the Riemannian metric $g$ of $M$ is fixed, reduction of holonomy
to the natural structure group is of limited interest as a criterion for
optimality in almost contact geometry.  A more promising approach is to
seek those almost contact structures whose characteristic vector field
$\xi$ is a harmonic section of the unit tangent bundle of $M$ \cite{28}. 
The study of such ``harmonic unit vector fields'' has been very active
over the past few years, as evidenced by the extensive bibliography of
\cite{11}.  In particular, contact metric manifolds whose characteristic
vector field is harmonic were recently characterized in \cite{22}, and
observed to constitute a very large class.  This reflects the incomplete
description of an almost contact structure by its characteristic vector
field; equally significant is the almost complex structure $J$ induced in
the almost contact subbundle $\D=\uperp\xi$ by the almost contact tensor
$\phi$.  In order to weigh the entire almost contact structure it is
natural to consider the corresponding section $\s$ of the associated
homogeneous fibre bundle $\pi\colon N\to M$ with fibre
$SO(2n+1)/U(n)$, which may be regarded as the odd-dimensional analogue of
the twistor bundle in almost Hermitian geometry.  We then look for those
$\s$ which are {\sl harmonic sections\/} of $\pi$.  This
methodology has already been applied to the study of almost complex
structures \cite{27, 8} and almost product structures \cite{26, 6, 77} (of
which harmonic unit vector fields are a special case), and recently a
general theory of ``harmonic reductions'' was laid out in \cite{29}. 
Following
\cite{29}, in \S2 we study the geometry of the {\sl universal almost
contact structure\/} for a Riemannian vector bundle, and in \S3 we use this
to obtain the harmonic section equations for $\s$, in terms of $\xi$ and
$J$.  The calculations are more complicated and subtle than those for
previously studied structures because the fibre of $\pi$ is no longer an
irreducible symmetric space; in fact it is neither irreducible nor
symmetric.  A naive expectation is that the harmonic section equations 
reduce to both $\xi$ and $J$ being harmonic.  This turns out to be false:
there is a first order coupling between $\xi$ and $J$ which in general
prevents such a reduction (Theorem 3.2).  We also derive the equations for
$\s$ to be a harmonic map (Theorem 3.4); as pointed out in \cite{29} in
general, these involve the curvature of $(M,g)$.
\par
In \S4 we examine the simplest situation where the harmonic section
equations for $\s$ reduce to the harmonicity of $\xi$.  This is when $J$ is
parallel, when we say that $\D$ is a {\sl K\"ahler bundle.} 
Examples include all $3$-dimensional almost contact structures, real
hypersurfaces of K\"ahler manifolds (Theorem 4.3), and all trans-Sasakian structures (Theorem 4.8).  The geometric data used for Theorem 4.8 yields an alternative proof of Marrero's
structure theorem \cite{19}, using Bianchi's first identity (Theorem 4.7). 
\par
To study examples where $\D$ is not a K\"ahler bundle, in \S5 we
assume that $M$ is an isometrically immersed hypersurface of an almost
Hermitian manifold $\tilde M$, and equipped with the induced almost contact
structure.  We relate the harmonic section equations for the almost contact
structure with those for the ambient almost Hermitian structure
(Proposition 5.3).  This allows us to study the canonical hypersurface
embedding, where $\tilde M=M\times\Bbb R$ (Theorem 5.4).  We use this to
show that any {\sl nearly cosymplectic\/} structure with parallel
characteristic field is parametrized by a harmonic map.  We also study the
case when $\tilde M$ is nearly K\"ahler.  Here we are able to characterize
the harmonicity of the induced almost contact structure when $M$ is a
contact metric manifold (Theorem 5.9), or a totally umbilic submanifold
(Theorem 5.10).  In particular, $\s$ is a harmonic map for the nearly
cosymplectic $S^5$ (Theorem 5.11).  We also exhibit examples where $\xi$ is
harmonic but
$\s$ is {\it not\/} a harmonic section; these include the {\sl nearly
Sasakian\/} $S^5$
\cite{5}.
\par
Further examples of harmonic almost contact structures are considered in
\cite{24}.
\bigskip
\head
2. Universal Almost Contact Geomtry
\endhead
For clarity we consider the general setup of an
oriented vector bundle $\E\to M$ of odd rank $r=2k+1$, with connection
$\nabla$ and holonomy-invariant fibre metric $\<\,,\>$.  Thus we no longer
require that $M$ is odd-dimensional, although our primary application will
be when this is the case and $\E=TM$.
\par
Let $G=SO(2k+1)$ and $H=U(k)$, included as follows:
$$
A+iB\mapsto
\pmatrix A&-B&0\\B&A&\vdots\\0&\cdots&1\endpmatrix
$$
Let $\phi_0$ be the skew-symmetric endomorphism with matrix:
$$
\pmatrix\Bbb O_k&-\Bbb I_k&0\\
\Bbb I_k&\Bbb O_k&\vdots\\0&\cdots&0
\endpmatrix
$$
Then $\phi_0$ is an element of the Lie algebra $\fg=\Skew(\Bbb R^r)$ of $G$,
and
$H$ and its Lie algebra
$\fh$ may be characterized:
$$
H\,=\,\{A\in G:A\,\phi_0\'A^{-1}=\phi_0\},
\qquad
\fh\,=\{a\in\fg:[\'a,\phi_0\']=0\},
$$
where $[\,,]$ is the commutator.  Let $\fm\subset\fg$ be the orthogonal
complement of $\fh$ with respect to the Killing form of $\fg$, and let
$\fm_1,\fm_2\subset\fm$ be the  following subspaces:
$$
\fm_1\,=\,\bigl\{a\in\fg:\{\'a,\phi_0\}=0\bigr\},\qquad
\fm_2\,=\,\bigl\{\{a,\eta_0\otimes\xi_0\}:a\in\fg\bigr\},
$$
where $\{\,,\}$ is the anticommutator, 
$\xi_0=(0,\dots,0,1)\in\Bbb R^r$, and $\eta_0\in(\Bbb R^r)^*$ is dual
to $\xi_0$.  Then $\fg=\fh\oplus\fm_1\oplus\fm_2$ is an
$\Ad(H)$-invariant splitting, with relations:
$$
\gather
[\'\fh,\fm_1\']\,\i\,\fm_1,\qquad
[\'\fh,\fm_2\']\,\i\,\fm_2 \\
[\'\fm_1,\fm_1\']\i\fh,\qquad
[\'\fm_2,\fm_2\']\,\i\,\fh\oplus\fm_1,\qquad
[\'\fm_1,\fm_2\']\i\fm_2
\tag2-1
\endgather
$$
If $a\in\fg$ then the decomposition
$a=a_\fh+a_{\fm_1}+a_{\fm_2}$ into components is:
$$
\gathered
a_\fh\,=\,-\tfrac12\bigl(\phi_0\{a,\phi_0\}\,
+\,a\circ(\eta_0\otimes\xi_0)\bigr), \\
\vspace{1ex}
a_{\fm_1}\,=\,\tfrac12\bigl(\phi_0\'[\'a,\phi_0\']\,
-\,a\circ(\eta_0\otimes\xi_0)\bigr),\qquad
a_{\fm_2}\,=\,\{a,\eta_0\otimes\xi_0\}
\endgathered
\tag2-2
$$
The following $H$-equivariant isomorphisms will also be useful:
$$
\aligned
\fh\oplus\fm_1&\to\Skew(\Bbb R^{2k}):
a\mapsto\hat a=a\'|\,\Bbb R^{2k} \\
\fm_2&\to\Bbb R^{2k}:
a_2\mapsto\hat a_2=a_2(\xi_0)
\endaligned
\tag2-3
$$
The first of these is an isometry, whereas the second is homothetic:
$$
|\'a_2\'|^2\,=\,2\'|\'\hat a_2\'|^2
\tag2-4
$$  
\par
Let $\mu\colon P\to M$ be the principal $G$-bundle of positively oriented
orthonormal frames of $\E$.  Then $\mu=\pi\circ\nu$ where
$\nu\colon P\to P/H=N$ is a principal $H$-bundle, and $\pi\colon N\to M$ is
naturally isomorphic to the homogeneous $G/H$-bundle associated to $\mu$. 
Let $\bE=\pi^*\E$ be the pullback bundle, equipped with the pullback
Riemannian structure.  Thus $\bE\to N$ is a vector bundle with reduced
structure group $H$.  Since $\xi_0$ (resp\. $\phi_0$) is $H$-invariant there
exists an induced unit section $\bxi$ of $\bE$ (resp\. skew-symmetric
endomorphism $\bPhi$ of $\bE$).  Then $\ker\bPhi=\<\bxi\>$, the line
subbundle generated by $\bxi$, and we define $\bD=\im\bPhi$.  The
restriction $\bPhi\'|\'\bD=\bJ$, a complex structure in $\bD$.
\par
Let $\bfH\to N$ and $\bfM_i\to N$ be the vector bundles associated to $\nu$
whose fibres are the $H$-modules $\fh$ and $\fm_i$.  Each is a subbundle of
$\Skew(\bE)$: 
\par\qquad
$\bfH$ (resp\. $\bfM_1$) commutes (resp\.
anticommutes) with $\bPhi$;
\par\qquad
$\bfM_2$ interchanges $\bD$ and $\<\bxi\>$.  
\flushpar
It follows that if $\a$ lies in $\bfH\oplus\bfM_1$ then $\a(\bxi)=0$.  The
$H$-equivariant ``hat'' isomorphisms (2-3) induce isomorphisms of $\bfH$
and $\bfM_i$ with vector bundles $\widehat\bfH\to N$ and $\widehat\bfM_i\to
N$ respectively, where $\widehat\bfH$ (resp\. $\widehat\bfM_1$) are
the subbundles of $\Skew(\bD)$ which commute (resp\. anticommute)
with $\bJ$, and $\widehat\bfM_2=\bD$.  If $\a$ lies in $\bfH\oplus\bfM_1$
then $\hat\a=\a\'|\'\bD$, whereas if $\a$ is in $\bfM_2$ then
$\hat\a=\a(\bxi)$.  For all sections $\a_i,\b_i$ of $\bfM_i$ it is easy to
show that:
$$
\align
\a_2(u)\,&=\,-\<\hat\a_2,u\>\bxi,
\quad\text{for all $u\in\bD$} 
\tag2-5 \\
[\'\bPhi,\a_2\']\hat{\phantom{2}}
&=\,\bJ(\hat\a_2)
\tag2-6 \\
[\'\a_1,\a_2\']\hat{\phantom{2}}
&=\,\hat\a_1(\hat\a_2)
\tag2-7 
\endalign
$$
It follows from (2-6) that:
$$
\text{\sl{if\/ $\a_2=[\'\bPhi,\b_2\']\;$ then $\;\b_2=-[\'\bPhi,\a_2\']$.}}
\tag2-8
$$
\par
If $\bfM=\bfM_1\oplus\bfM_2\to N$ then the {\sl homogeneous connection
form\/} (cf\. \cite{29}) is the $\bfM$-valued 1-form $\theta$ on $N$
obtained by projecting the $\fm$-component of the connection form
$\omega\in\Omega^1(P,\fg)$.  We note that $\nabla\bPhi$ is also an
$\bfM$-valued $1$-form on $N$.  For, if $\boldeta$ is the section of $\bE^*$
dual to $\bxi$ then differentiating the identity 
$$
\bPhi^2\,=\,-1\,+\,\boldeta\otimes\bxi
\tag2-9
$$ 
yields:
$$
\{\'\nabla\bPhi,\bPhi\'\}\,
=\,\nabla\boldeta\otimes\bxi\,+\,\boldeta\otimes\nabla\bxi
$$
Hence:
$$
\align
\bPhi\'\{\nabla\bPhi,\bPhi\}\,
+\,\nabla\bPhi\circ(\boldeta\otimes\bxi)\,
&=\,\boldeta\otimes\left(\bPhi\,\nabla\bxi+\nabla\bPhi(\bxi)\right)\,
=\,\boldeta\otimes\nabla(\bPhi\'\bxi)\,=\,0,
\endalign
$$
and it follows from (2-2) that the $\bfH$-component of
$\nabla\bPhi$ vanishes.  Both $\theta$ and
$\nabla\bPhi$ therefore split into $\bfM_i$-components: 
$$
\theta\,=\,\theta_1+\theta_2,\qquad
\nabla\bPhi\,=\,(\nabla\bPhi)_1\,+\,(\nabla\bPhi)_2
$$
\proclaim{2.1 Lemma}
\flushpar
{\rm(a)}\quad
$\theta_1=\frac12\'\bPhi\circ(\nabla\bPhi)_1$\rom;
\qquad
{\rm(b)}\quad
$\theta_2=\left[\'\bPhi,(\nabla\bPhi)_2\right]$
\endproclaim
\demo{Proof}
Since $\bPhi$ lifts to the (constant) $H$-equivariant map
$\phi_0\colon P\to\fh$ we have (cf\. \cite{29, Proposition 4.1}):
$$
\nabla\bPhi\,=\,[\'\theta,\bPhi\']\,
=\,[\'\theta_1,\bPhi\']\,+\,[\'\theta_2,\bPhi\']\,
=\,-2\'\bPhi\circ\theta_1\,+\,[\'\theta_2,\bPhi\'],
$$
and therefore:
$$
(\nabla\bPhi)_1\,=\,-2\'\bPhi\circ\theta_1,
\qquad
(\nabla\bPhi)_2\,=\,[\'\theta_2,\bPhi\'],
$$
from which (a) follows directly, and (b) follows via (2-8).
\qed
\enddemo
An expression for $\theta$ in terms of $\bJ$ and $\bxi$ is obtained by
applying the ``hat'' isomorphisms to Lemma 2.1.  Let $\bar\nabla$
denote the connection in $\bD\to N$ obtained by orthogonally projecting
$\nabla$, and note that $\nabla\bxi$ is $\bD$-valued  since $\bxi$ is a
unit section.
\proclaim{2.2 Proposition}
\flushpar
{\rm(a)}\quad
$\hat\theta_1\,=\,\frac12\,\bJ\,\bar\nabla\bJ$\rom;
\qquad
{\rm(b)}\quad
$\hat\theta_2\,=\,\nabla\bxi$
\endproclaim
\demo{Proof}
\flushpar
(a)\quad
If $A\in TN$ and $u$ is a section of $\bD$ then since $\bfM_2$
swaps $\bD$ and $\<\bxi\>$ it follows from Lemma 2.1\,(a) that:
$$
\align
2\hat\theta_1(A)(u)\,
&=\,2\theta_1(A)(u)\,
=\,\bPhi\'(\nab A\bPhi)_1(u)\,
=\,\bPhi\,\nab A\bPhi(u) \\
&=\,\bPhi\,\nab A(\bJ u)\,+\,\nab A u\,-\,\<\nab A u,\bxi\>\bxi,
\quad\text{by (2-9)} \\
&=\,\bJ\,\bar\nab A(\bJ u)\,+\,\bar\nab A u\,
=\,\bJ(\bar\nab A(\bJ u)\,-\,\bJ\,\bar\nab A u)\,
=\,\bJ\,\bar\nab A\bJ(u).
\endalign
$$
(b)\quad
By Lemma 2.1\,(b) and (2-6):
$$
\align
\hat\theta_2\,&=\,\bJ\'(\nabla\bPhi)_2\hat{\vphantom{1}}\,
=\,\bPhi\'(\nabla\bPhi)_2(\bxi)\,
=\,\bPhi\left\{\nabla\bPhi,\boldeta\otimes\bxi\right\}(\bxi),
\quad\text{by (2-2)} \\
&=\,\bPhi\left(\nabla\bPhi(\bxi)\,
+\,\left\<\nabla\bPhi(\bxi),\bxi\right\>\bxi\right)\,
=\,-\bPhi^2\,\nabla\bxi\,=\,\nabla\bxi,
\quad\text{by (2-9).}
\tag"$\ssize\square$"
\endalign
$$
\enddemo
The $\fh$-component of $\omega$ is a connection form in $\nu\colon P\to N$,
and induces linear connections in the associated vector bundles 
$\bfM,\bfM_i\to N$, each of which we denote by $\nabla^c$ and call the {\sl
canonical connection\/} for the bundle in question.  For all sections
$\a_i$ of $\bfM_i$:
$$
\nabla^c(\a_1+\a_2)\,=\,\nabla^c\a_1\,+\,\nabla^c\a_2
\tag2-10
$$
\proclaim{2.3 Proposition}
For all sections $\a_i$ of $\bfM_i$ we have:
\flushpar
{\rm(a)}\quad
$(\nabla^c\a_1)\hat{\phantom{2}}
=\,\frac12\'\bJ\,[\'\bar\nabla\hat\a_1,\bJ\']$\rom;
\qquad
{\rm(b)}\quad
$(\nabla^c\a_2)\hat{\phantom{2}}
=\,\bar\nabla\hat\a_2\,
-\,\frac12\'\bJ\,\bar\nabla\bJ(\hat\a_2)$
\endproclaim
\demo{Proof}
If $\a$ is a section of $\bfM$ then by the structure equations
(cf\. \cite{29, Proposition 2.7}):
$$
\nabla\a\,=\,\nabla^c\a\,+\,[\'\theta,\a\']\,
=\,\nabla^c\a\,+\,[\'\theta_1,\a\']\,+\,[\'\theta_2,\a\']
\tag2-11
$$
\flushpar
(a)\quad
By (2-11) and (2-1), $\nabla^c\a_1$ is
the $\bfM_1$-component of $\nabla\a_1$,  and hence by (2-2):
$$
2\,\nabla^c\a_1\,=\,\bPhi\'[\'\nabla\a_1,\bPhi\']\,
-\,\nabla\a_1\circ(\boldeta\otimes\bxi)
$$
Therefore if $u\in\bD$ then:
$$
\align
2(\nabla^c\a_1)\hat{\phantom{2}}(u)\,
&=\,2\,\nabla^c\a_1(u)\,
=\,\bPhi\'[\'\nabla\a_1,\bPhi\'](u)\,
=\,\bPhi\,\nabla\a_1(\bPhi u)\,
-\,\bPhi^2\,\nabla\a_1(u) \\
&=\,\bPhi\bigl(\bar\nabla\hat\a_1(\bJ u)\,
+\,\<\nabla\a_1(\bJ u),\bxi\>\'\bxi\bigr)\,
+\,\nabla\a_1(u)\,-\,\<\nabla\a_1(u),\bxi\>\'\bxi \\
&=\,\bJ\,\bar\nabla\hat\a_1(\bJ u)\,-\,\bar\nabla\hat\a_1(u)\,
=\,\bJ\,[\'\bar\nabla\hat\a_1,\bJ\'](u)
\endalign
$$
(b)\quad
It follows from (2-11), (2-1) and (2-2) that:
$$
\nabla^c\a_2\,+\,[\'\theta_1,\a_2\']\,
=\,\left\{\nabla\a_2,\boldeta\otimes\bxi\right\}
$$
Each term in this equation is $\bfM_2$-valued, so for all
$A$ in $TN$ we have:
$$
\align
(\nabc A\a_2)\hat{\phantom{2}}\,
&=\,\{\nab A\a_2,\boldeta\otimes\bxi\}(\bxi)\,
-\,\hat\theta_1(A)(\hat\a_2),
\quad\text{by (2-7)} \\
&=\,\nab A\a_2(\bxi)\,
-\,\tfrac12\,\bJ\,\bar\nab A\bJ(\hat\a_2),
\quad\text{by Proposition 2.2\,(a),} 
\endalign
$$
since $\nab A\a_2$ is skew-symmetric.  Now note that:
$$
\align
\nabla\a_2(\bxi)\,
&=\,\nabla\hat\a_2\,-\,\a_2(\nabla\bxi)\,
=\,\nabla\hat\a_2\,+\,\<\hat\a_2,\nabla\bxi\>\bxi,
\quad\text{by (2-5)} \\
&=\,\nabla\hat\a_2\,-\,\<\nabla\hat\a_2,\bxi\>\bxi\,
=\,\bar\nabla\hat\a_2
\tag"$\ssize\square$"
\endalign
$$
\enddemo
Finally, we recall from \cite{29} that the {\sl homogeneous curvature
form\/} is the $\bfM$-valued $2$-form $\Theta$ on $N$ obtained by projecting
the $\fm$-component $\Omega_\fm$ of the curvature form
$\Omega\in\Omega^2(P,\fg)$.  Then $\Theta=\Theta_1+\Theta_2$. 
Let $R$ be the curvature tensor of $\nabla$.  Since $\Theta_i$ (resp\. $R$)
is the projection to $N$ (resp\. $M$) of $\Omega_{\fm_i}$ (resp\. $2\Omega$)
it follows from (2-2) that:
$$
4\Theta_1\,=\,\bPhi\'[\'\pi^*R,\bPhi\']\,
-\,\pi^*R\circ(\boldeta\otimes\bxi),\qquad
2\Theta_2\,=\,\{\'\pi^*R,\boldeta\otimes\bxi\'\}
\tag2-12
$$
\bigskip
\head
3. The Harmonic Section Equations
\endhead
Let $\s\colon M\to N$ be a section of $\pi$.  The corresponding
almost contact structure in $\E$ is obtained by pulling back the universal
almost contact structure:
$$
\phi\,=\,\s^*\bPhi,\qquad
\xi\,=\,\s^*\bxi
$$
Then $\D=\s^*\bD$ and $J=\s^*\bJ$.
Let $Q\to M$ be the principal $H$-subbundle of $P\to M$ comprising all
{\sl $\phi$-adapted\/} positively oriented orthonormal frames of $\E$, and
let $\fM\to M$ (resp\. $\fM_i\to M$) be the associated vector bundle
with fibre $\fm$ (resp\. $\fm_i$).  Then:
$$
\fM\,=\,\s^*\bfM,\qquad
\fM_i\,=\,\s^*\bfM_i
$$
Thus $\fM_1$ (resp\. $\fM_2$) is the subbundle of $\Skew(\E)$ which
anticommutes with $\phi$ (resp\. swaps $\D$ and $\<\xi\>$).  Applying the
``hat'' isomorphisms yields the bundle $\widehat\fM_1$ of skew-symmetric
endomorphisms of $\D$ which anticommute with $J$, and $\widehat\fM_2=\D$.
Define $\psi=\s^*\theta$, an $\fM$-valued 1-form on $M$.  Pulling back
Proposition 2.2 by $\s$ yields:
$$
\hat\psi_1\,=\,\tfrac12\'J\,\bar\nabla J,
\qquad
\hat\psi_2\,=\,\nabla\xi
\tag3-1
$$
In particular, $\s$ is horizontal ($\psi=0$) if and only if $\xi$ and $J$
are parallel.
\par
Pulling back $\nabla^c$ by $\s$  yields connections in $\fM,\fM_i\to M$, which will also be
denoted by $\nabla^c$; they are simply the linear connections in these
bundles associated to the $H$-connection form $\omega_\fh\'|\'_{TQ}$.  
For all sections $\a_i$ of $\fM_i$ it follows from Proposition 2.3 that:
$$
(\nabla^c\a_1)\hat{\vphantom{1}}\,
=\,\tfrac12\'J\,[\'\bar\nabla\hat\a_1,J\'],
\qquad
(\nabla^c\a_2)\hat{\vphantom{1}}\,
=\,\bar\nabla\hat\a_2-\tfrac12\'J\,\bar\nabla J(\hat\a_2)
\tag3-2
$$ 
We note that $\,\nabla^c\psi=\nabla^c\psi_1+\nabla^c\psi_2$ by (2-10), and introduce
$2$-covariant tensors on $M$:
$$
\<\nabla\xi\otimes\nabla\xi\>(X,Y)\,
=\,\<\nab X\xi,\nab Y\xi\>,
\qquad
(\bar\nabla J\otimes\nabla\xi)(X,Y)\,
=\,\bar\nab XJ(\nab Y\xi),
$$
the first (resp\. second) of which is symmetric (resp\. $\D$-valued).
\proclaim{3.1 Proposition}
\flushpar
{\rm(a)}\quad
$(\nabla^c\psi_1)\hat{\vphantom{1}}\,
=\,\frac14\'[\'J,\bar\nabla^2J\']$
\flushpar
{\rm(b)}\quad
$(\nabla^c\psi_2)\hat{\vphantom{1}}\,
=\,\nabla^2\xi\,
+\,\<\nabla\xi\otimes\nabla\xi\>\xi\,
-\,\frac12\'J\,\bar\nabla J\otimes\nabla\xi$
\endproclaim
\demo{Proof}
\flushpar
(a)\quad
From the first of equations (3-2):
$$
\align
4(\nabla^c\psi_1)\hat{\vphantom{1}}\,
&=\,2\,J\'\bigl[\'\bar\nabla\hat\psi_1,J\'\bigr]\,
=\,J\,\bigl[\bar\nabla(J\,\bar\nabla J),J\'\bigr],
\quad\text{by (3-1)} \\
&=\,J\,[\,\bar\nabla J\circ\bar\nabla J,J\']\,
+\,J\,[\'J\,\bar\nabla^2J,J\']\,=\,J\,[\'J\,\bar\nabla^2J,J\'],
\quad\text{since $\{\bar\nabla J,J\}=0$,} \\
&=\,J^2\circ\bar\nabla^2J\circ J\,-\,J^3\circ\bar\nabla^2J\,
=\,J\circ\bar\nabla^2J\,-\,\bar\nabla^2J\circ J\,
=\,[\'J,\bar\nabla^2J\'].
\endalign
$$
(b)\quad
From the second of equations (3-2), and (3-1):
$$
(\nabla^c\psi_2)\hat{\vphantom{1}}\,
=\,\bar\nabla\hat\psi_2\,-\,\tfrac12\'J\,\bar\nabla J(\hat\psi_2)\,
=\,\bar\nabla(\nabla\xi)\,-\,\tfrac12\'J\,\bar\nabla J\otimes\nabla\xi,
$$
and since $\xi$ is a unit vector:
$$
\align
\bar\nab X(\nabla\xi)(Y)\,
&=\,\nab X(\nabla\xi)(Y)\,-\,\bigl\<\nab X(\nabla\xi)(Y),\xi\bigr\>\xi \\
&=\,\nabsq XY\xi\,-\,\bigl\<\nab X\nab Y\xi,\xi\bigr\>\xi\,
=\,\nabsq XY\xi\,+\,\bigl\<\nab X\xi,\nab Y\xi\bigr\>\xi.
\tag"\qed"
\endalign
$$
\enddemo
Let $\V\to N$ be the $\pi$-vertical subbundle of $TN$.  
The {\sl vertical tension field\/} of $\s$ is:  
$$
\tau^v\s\,=\,\Tr\nabla^v\dvs\,
=\,\nabv{E_i}\dvs(E_i),
$$
where $\{E_i\}$ is any orthonormal frame of $TM$ (here and henceforward
the summation convection is in force), and
$\dvs$ (resp\. $\nabla^v$) is the $\V$-component of $d\s$ (resp\. the
Levi-Civita connection of $N$).  (The Riemannian metric on $N$ is obtained
by horizontally lifting $g$, and supplementing with the fibre metric in
$\V$ induced by the $G$-invariant metric on $G/H$.)  Then $\s$ is a
{\sl harmonic section\/} of $\pi$ if and only if $\tau^v\s=0$.  From
\cite{29}, the restriction $\theta\'|\'_\V$ is a canonical vector bundle
isomorphism $I\colon\V\to\bfM$, and $I\circ\dvs=\psi$.  Because $G/H$ is not
a symmetric space, $I$ is not connection-preserving;  however since $G/H$
is naturally reductive we have:
$$
I\circ\tau^v\s\,=\,-\delta^c\psi\,=\,\nabc{E_i}\psi(E_i),
$$
by \cite{29, Theorem 3.2}.  We have $I=I_1+I_2$ from the splitting of 
$\bfM$, and by (2-10):
$$
I_i\circ\tau^v\s\,=\,-\delta^c\psi_i
$$
We also recall the definitions of the {\sl rough Laplacians:}
$$
\bar\nabla^*\bar\nabla J\,=\,-\bar\nabsq{E_i}{E_i}J,
\qquad
\nabla^*\nabla\xi\,=\,-\nabsq{E_i}{E_i}\xi
$$
In general,  a Hermitian structure $J$ in a Riemannian vector bundle is said
to be {\sl harmonic\/} if the corresponding section of the associated
twistor bundle is a harmonic section, and by \cite{27} this is the case
precisely when $J$ commutes with its rough Laplacian.  We also say that
{\sl $\xi$ is harmonic\/} if $\xi$ is a harmonic section of
the unit sphere bundle  $U\E\to M$, and by \cite{28} this is the case
precisely when:
$$
\nabla^*\nabla\xi\,=\,|\'\nabla\xi\'|^2\xi
\tag3-3
$$
We define a section $T(\phi)$ of $\D$ as follows:
$$
T(\phi)\,=\,\Tr(\bar\nabla J\otimes\nabla\xi)\,
=\,\bar\nab{E_i}J(\nab{E_i}\xi)
$$
and abbreviate:
$$
\bar\tau(J)\,=\,\tfrac14\,[\'\bar\nabla^*\bar\nabla J,J\']
$$
Thus $\bar\tau(J)$ is a skew-symmetric endomorphism of $\D$ which
anticommutes with $J$. The following result is an immediate consequence of
Proposition 3.1.  
\proclaim{3.2 Theorem}
\flushpar
{\rm(a)}\quad
$(I_1\,\tau^v\s)\hat{\vphantom{1}}\,
=\,\bar\tau(J)$
\smallskip\noindent
{\rm(b)}\quad
$(I_2\,\tau^v\s)\hat{\vphantom{1}}\,
=\,-\nabla^*\nabla\xi\,+\,|\'\nabla\xi\'|^2\'\xi\,
-\,\frac12\'J\,T(\phi)$
\flushpar
Therefore $\s$ is a harmonic section if and only if $J$ is harmonic and: 
$$
\nabla^*\nabla\xi\,=\,|\'\nabla\xi\'|^2\xi\,-\,\tfrac12\'J\,T(\phi)
$$
If $\xi$ and $J$ are harmonic then $\s$ is a harmonic section
precisely when $T(\phi)=0$.
\endproclaim
We now define an $\fM$-valued $2$-form $\Psi=\s^*\Theta$ on $M$.   Let
$R_\D$ denote the curvature-type tensor in $\D$ obtained by orthogonally
projecting $R$:
$$
R_\D(X,Y)u\,=\,R(X,Y)u\,-\,\<R(X,Y)u,\xi\>\xi,
\quad\text{for all $u\in\D$.}
$$
If $\bar R$ is the curvature tensor of the connection $\bar\nabla$ in $\D$
then:
$$
\bar R\,=\,R_\D\,+\,r(\nabla\xi,\nabla\xi),
\tag3-4
$$
where $r$ is the following curvature-type operator in $\E$:
$$
r(u,v)w\,=\,\<v,w\>u\,-\,\<u,w\>v,
\qquad\text{for all $u,v,w\in\E$.}
\tag3-5
$$
\proclaim{3.3 Proposition}
\flushpar
\rom{(a)}\quad
$\hat\Psi_1\,=\,\tfrac14J\'[\'R_\D,J\']$\rom;
\qquad
\rom{(b)}\quad
$\hat\Psi_2(X,Y)=\tfrac12R(X,Y)\xi$
\endproclaim
\demo{Proof}
\;(a)\quad
Pulling back equation the first of equations (2-12) by $\s$ yields:
$$
\align
4\hat\Psi_1(X,Y)u\,
&=\,4\Psi_1(X,Y)u\,
=\,\phi\'[\'R(X,Y),\phi\']u\,
=\,J\'[R_\D(X,Y),J\']u
\endalign
$$
(b)\quad
From the second of equations (2-12):
$$
\align
2\hat\Psi_2(X,Y)\,
&=\,2\Psi_2(X,Y)\xi\,
=\,\{R(X,Y),\eta\otimes\xi\}\xi\,
=\,R(X,Y)\xi.
\tag"\qed"
\endalign
$$
\enddemo 
We recall from \cite{29} that $\s$ is a {\sl harmonic map\/}
precisely when $\s$ is a harmonic section and the following $1$-form
$\<\psi,\Psi\>$ on $M$ vanishes:
$$
\align
\<\psi,\Psi\>(X)\,&=\,\<\psi\otimes\Psi\>(E_i,E_i,X)\,
=\,\<\psi(E_i),\Psi(E_i,X)\> \\
&=\,\<\psi_1,\Psi_1\>(X)\,+\,\<\psi_2,\Psi_2\>(X)
\endalign
$$
If $\bar\s$ is the section of the twistor bundle of $\D$ parametrizing $J$
then from \cite{29}:
$$
\<\bar\psi,\bar\Psi\>(X)\,
=\,\<\bar R(E_i,X),J\,\bar\nab{E_i}J\>,
\tag3-6
$$
and if $\xi$ is harmonic then $\xi\colon M\to U\E$ is a harmonic
map precisely when \cite{10}:
$$
\<R(E_i,X)\xi,\nab{E_i}\xi\>\,=\,0
\tag3-7
$$
Use of Proposition 3.3 with (3-2), (3-4) and (2-4) yields the following:
\proclaim{3.4 Theorem}
For all $X$ in $TM$ we have:
\TagsOnLeft
$$
\align
\<\psi_1,\Psi_1\>(X)\,&=\,\tfrac14\<R(E_i,X),J\,\bar\nab{E_i}J\>\,
=\,\<\bar\psi,\bar\Psi\>(X)\,
+\,\tfrac12\<JT(\phi),\nab X\xi\> 
\tag a \\
&=\,\<\bar\psi,\bar\Psi\>(X)\,-\,\<\nabla^*\nabla\xi,\nab X\xi\>,
\quad\text{if $\s$ is a harmonic section.}  \\
\vspace{1ex}
\<\psi_2,\Psi_2\>(X)\,
&=\,\<R(E_i,X)\xi,\nab{E_i}\xi\>
\tag b
\endalign
$$
\TagsOnRight
If $\s$ is a harmonic section, and $\bar\s$ and $\xi$ are
harmonic maps, then $\s$ is a harmonic map.
\endproclaim
\bigskip
\head 
4. K\"ahler Contact Bundle
\endhead 
There are two partial reductions of the holonomy of $\E$,
which yield significant simplifications to Theorems 3.2 and 3.4.
\definition{Definition}
A Riemannian vector bundle equipped with a parallel Hermitian structure is
said to be a {\sl K\"ahler bundle.}
\enddefinition
\proclaim{4.1 Theorem}
\flushpar
\rom{(a)}\quad
If $(\D,\bar\nabla,J)$ is a K\"ahler bundle then $\s$ is a harmonic section
\rom(resp\. harmonic map\rom) precisely when $\xi$ is a harmonic section of
$U\E$ \rom(resp\. harmonic map\rom).
\flushpar
\rom{(b)}\quad
If $\xi$ is parallel then $\s$ is a harmonic section \rom(resp\. harmonic
map\rom) precisely when $\bar\s$ is a harmonic section \rom(resp\. harmonic
map\rom).
\endproclaim
Henceforward we assume that $M$ has odd dimension $2n+1$, and $\E=TM$,
with the Levi Civita connection.  There are two ``obvious'' cases
when $\D$ is a K\"ahler bundle.  If $M$ is an oriented
(real) hypersurface of an almost Hermitian manifold $(\tilde M,g,J)$, with
compatible unit normal $\nu$, then there is an induced almost contact
structure:
$$
\xi\,=\,-J\nu,
\qquad
\phi X\,=\,JX\,-\,\eta(X)\nu
\tag4-1
$$
The induced almost complex structure in $\D$ is thus the
restriction of $J$.
\proclaim{4.2 Lemma}
The induced almost complex structure in $\D$ satisfies:
$$
\bar\nabla J\,
=\,\tilde\nabla J\,
-\,\<\tilde\nabla J,\xi\>\xi\,
-\,\<\tilde\nabla J,\nu\>\nu
$$
\endproclaim
\demo{Proof}
This follows from the simple identity, for all sections $Z$ of $\D$:
$$
\bar\nabla Z\,=\,\tilde\nabla Z\,-\,\<\tilde\nabla Z,\xi\>\xi\,
-\,\<\tilde\nabla Z,\nu\>\nu
\tag"\qed"
$$
\enddemo
\proclaim{4.3 Theorem}
Suppose either of the following holds:
\flushpar
\rom{(a)}\quad
$M$ is a $3$-dimensional almost contact manifold. 
\flushpar
\rom{(b)}\quad
$M$ is an oriented real hypersurface of a K\"ahler manifold, equipped with
the induced almost contact structure.
\flushpar
Then $\D$ is a K\"ahler bundle.  Hence $\s$ is a harmonic section
\rom(resp\. harmonic map\rom) if and only if $\xi$ is a harmonic unit
vector field
\rom(resp\. harmonic map\rom).
\endproclaim
\par\pagebreak
\demo{Proof}
Part (b) follows directly from Lemma 4.2.  For (a) note that
if $Z\in\D_x$ is a unit vector, then $(Z,JZ)$ is an orthonormal basis of
$\D_x$, and for all $X\in\D_x$ we have:
$$
\<\bar\nab XJ(Z),Z\>\,=\,0,
\quad
2\<\bar\nab XJ(Z),JZ\>\,
=\,\<\bar\nab XJ(Z),JZ\>-\<JZ,\bar\nab XJ(Z)\>\,=\,0.
\tag"\qed"
$$
\enddemo
\example{4.4 Example}
Harmonic unit vector fields on $3$-manifolds have been extensively studied,
notably in \cite{12}; in  particular, each Thurston geometry has a natural
harmonic unit field.  Those on $S^2\times\Bbb R$, $H^2\times\Bbb R$, $\Nil$,
$\Sol$ and $\widetilde{SL_2}$ are invariant under the full isometry group. 
Every discrete subgroup of Euclidean isometries acting freely on $\Bbb R^3$
leaves a direction invariant, and the Hopf field on $S^3$ is invariant
under all finite subgroups of isometries which act freely \cite{23}.  Hence
by Theorem 4.3\,(a), any 3-manifold with non-hyperbolic geometric structure
has a natural harmonic almost contact structure.
\endexample
\example{4.5 Example}
Let $M=S^{2n+1}$, included in $\Bbb R^{2n+2}\cong\Bbb C^{n+1}$ as the unit
sphere, and oriented by the unit outward-pointing normal.  Then $-\xi$ is
the standard Hopf vector field, and the induced almost contact structure is
the standard Sasakian structure on $S^{2n+1}$. It is well-known that $\xi$
is a harmonic unit vector field \cite{25}, and indeed a harmonic map
\cite{16}; therefore $\s$ is a harmonic map by Theorem 4.3\,(b). 
\endexample
Our next example generalizes Example 4.5 in a different way.  It follows
from
$$
\bar\nab XJ(Z)\,
=\,\nab X\phi(Z)\,-\,\<\nab X\phi(Z),\xi\>\xi,
\qquad
Z\in\D,
\tag4-2
$$
that $\D$ is a K\"ahler bundle if and only if $\nab X\phi\in\fM_2$ for all
$X$ in $TM$.  This is clearly the case for a {\sl trans-Sasakian
structure of type $(\a,\b)$\/} \cite{21,\,4}:
$$
\nab X\phi(Y)\,=\,\a\'r(\xi,X)Y\,+\,\b\'r(\xi,\phi X)Y,
$$
where $\a,\b\colon M\to\Bbb R$ are smooth functions, and $r$ is defined
in (3-5).  If $\a=0$ (resp\. $\b=0$) the structure is said to be {\sl
$\b$-Kenmotsu\/} (resp. {\sl $\a$-Sasakian\/}).  The characteristic vector
field of a trans-Sasakian manifold satisfies \cite{9,\,17}:
$$
\align
\nab X\xi\,&=\,-\a\phi X\,-\,\b\phi^2X
\tag4-3 \\
\nabsq XY\xi\,
&=\,(\a^2-\b^2)\eta(Y)X\,+\,2\a\b\'\eta(Y)\phi X\,
-\,(X.\a)\phi Y\,-\,(X.\b)\phi^2Y
\tag4-4 \\
R(X,Y)\xi\,
&=\,\phi\'r(X,Y)\nabla\a\,+\,\phi^2r(X,Y)\nabla\b\,
+\,(\a^2-\b^2)r(X,Y)\xi\,+\,2\a\b\'\phi\'r(X,Y)\xi
\tag4-5 \\
\Ric(\xi)\,
&=\,\phi\,\nabla\a\,-\,\phi^2\,\nabla\b\,
-\,2n\,\nabla\b\,+\,2n(\a^2-\b^2)\xi,
\tag4-6
\endalign
$$
where $\nabla\a,\nabla\b$ are the gradient vectors.   Of particular interest
is the $2$-form $\R(\Phi)$, where
$\R\colon\Omega^2(M)\to\Omega^2(M)$ is the {\sl curvature operator\/} and
$\Phi$ is the {\sl fundamental $2$-form:} 
$$
\Phi(X,Y)\,=\,\<X,\phi Y\>
$$
If $\{F_i\}$ is an orthonormal frame of $\D$, extended to a
frame $\{E_i\}$ of $TM$ with $E_{2n+1}=\xi$, then:
$$
\align
\R(\Phi)(X,Y)\,
&=\,-\tfrac12\<R(X,Y)E_i,\phi E_i\>\,
=\,-\tfrac12\<R(X,Y)F_i,\phi F_i\> 
\tag4-7 \\
&=\,\<R(X,F_i)\phi F_i,Y\>,
\tag4-8
\endalign
$$ 
where (4-8) follows from (4-7) by the symmetries of $R$,
including Bianchi's first identity.  Let $\R(\Phi)^\sharp$ denote the
corresponding endomorphism field.  Using (4-5) in each of (4-7) and (4-8)
yields two expressions for $\R(\Phi)$, and an alternative proof of
\cite{19, Theorem 4.1}.
\proclaim{4.6 Lemma}
On a trans-Sasakian manifold we have:
\flushpar
\rom{(a)}\quad
$\R(\Phi)^\sharp(\xi)\,
=\,\phi^2\,\nabla\a\,-\,\phi\,\nabla\b$
\flushpar
\rom{(b)}\quad
$\R(\Phi)^\sharp(\xi)\,
=\,-\phi^2\,\nabla\a\,-\,\phi\,\nabla\b\,-\,2n\,\nabla\a\,-\,4n\a\b\'\xi$
\endproclaim
\proclaim{4.7 Theorem}
For all trans-Sasakian structures:\quad$d\a(\xi)=-2\a\b$.
\flushpar
If $\dim M\geqs5$ then $\a\b=0$ and $\nabla\a=0$.  Thus if $M$
is connected it is either $\a$-Sasakian with $\a$ constant,
or $\b$-Kenmotsu.
\endproclaim
\demo{Proof}
Since $\R(\Phi)$ is skew-symmetric it follows from Lemma 4.6\,(b) that:
$$
0\,=\,\R(\Phi)(\xi,\xi)\,=\,-2n\<\nabla\a,\xi\>\,-\,4n\a\b,
$$
and therefore:
$$
\<\nabla\a,\xi\>\,=\,-2\a\b
\tag4-9
$$
Now, equating the two identities of Lemma 4.6 yields:
$$
\align
0\,&=\,\phi^2\,\nabla\a\,+\,n\,\nabla\a\,+\,2n\a\b\'\xi\,
=\,\phi^2\,\nabla\a\,-\,n\'\phi^2\,\nabla\a\,
+\,n\<\nabla\a,\xi\>\xi\,+\,2n\a\b\'\xi \\ &=\,(1-n)\phi^2\,\nabla\a\,
+\,2(n-1)\a\b\'\xi,
\quad\text{using (4-9).}
\endalign
$$
Therefore if $n>1$ then the $\D$-component of $\nabla\a$ vanishes
and $\a\b=0$.  It then follows from (4-9) that  $\nabla\a=0$.
\qed
\enddemo
\remark{Remark}
Identity (4-9) was noted in \cite{9}.
\endremark
\definition{Definition}
A vector field (resp\. $1$-form) on $M$ is {\sl characteristic\/} if it is
proportional to $\xi$ (resp\. $\eta$).
\enddefinition
\proclaim{4.8 Theorem}
Suppose $M$ is trans-Sasakian of type $(\a,\b)$.  Then $\s$ is a harmonic
section if and only if:
$$
\nabla\a\,=\,-\phi\,\nabla\b\,-\,2\a\b\'\xi,
$$
and $\s$ is a harmonic map if and only if in addition:
$$
\xi.(\b^2)\,=\,2\b(\a^2-\b^2)
$$
\rom{(a)}\quad
If $M$ is $\a$-Sasakian then $\s$ is a harmonic section if and only if $\a$
is locally constant, in which case $\s$ is a harmonic map.
\flushpar
\rom{(b)}\quad
If $M$ is $\b$-Kenmotsu, or $\dim M\geqs5$, then $\s$ is a harmonic section
if and only if $\nabla\b$ is characteristic,  and $\s$ is a harmonic map if
and only if $\,\nabla\b=-\b^2\xi$.  If $M$ is compact $\b$-Kenmotsu then
$\s$ is a harmonic map precisely when $M$ is cosymplectic. 
\endproclaim
\demo{Proof}
From (4-4):
$$
\nabla^*\nabla\xi\,
=\,\phi\,\nabla\a\,+\,\phi^2\,\nabla\b\,
+\,2n(\a^2+\b^2)\xi
$$
Therefore by (3-3), $\xi$ is harmonic if and only if:
$$
\phi\,\nabla\a\,+\,\phi^2\,\nabla\b\,=\,0,
$$
which is equivalent to:
$$
\align
0\,&=\,\phi^2\,\nabla\a\,-\,\phi\,\nabla\b\,
=\,-\nabla\a\,+\,\<\nabla\a,\xi\>\xi\,-\,\phi\,\nabla\b \\
&=\,-\nabla\a\,-\,2\a\b\'\xi\,-\,\phi\,\nabla\b,
\quad\text{by Theorem 4.7.}
\tag4-10
\endalign
$$
An $\a$-Sasakian (resp\. $\b$-Kenmotsu) structure satisfies (4-10) if and
only if $\a$ is locally constant (resp\. $\nabla\b$ is characteristic).  
If $\dim M\geqs5$ then (4-10) holds precisely when $\nabla\b$ is
characteristic, by Theorem 4.7.
\par
From Theorem 3.4 and (4-3):
$$
\align
\<\psi,\Psi\>(X)\,
&=\,\<R(E_i,X)\xi,\nab{E_i}\xi\>\,
=\,-\a\<R(F_i,X)\xi,\phi F_i\>\,
+\,\b\<R(F_i,X)\xi,F_i\> \\
&=\,\a\R(\Phi)(\xi,X)\,+\,\b\Ric(\xi,X) \\
&=\,\a\<\phi^2\,\nabla\a-\phi\,\nabla\b,X\>\,
+\,\b\<\phi\,\nabla\a\,-\,\phi^2\,\nabla\b\,
-\,2n\,\nabla\b\,+\,2n(\a^2-\b^2)\xi,X\>,
\endalign
$$
by (4-6) and Lemma 4.6\,(a).  Inserting the harmonic section equations
(4-10), and noting that when $n>1$ these imply that $\nabla\b$ is
characteristic:
$$
\align
\<\psi,\Psi\>^\sharp\,
&=\,2\b\big(n(\a^2-\b^2)\xi\,
-\,\phi^2\,\nabla\b\,-\,n\,\nabla\b\big) \\
&=\,2\b\big(n(\a^2-\b^2)\xi\,
+\,(n-1)\phi^2\,\nabla\b\,
-\,n\<\nabla\b,\xi\>\xi\big) \\
&=\,2n\b\big((\a^2-\b^2)\,-\,(\xi.\b)\big)\xi,
\quad\text{for all $n$.}
\endalign
$$
Therefore if $\s$ is a harmonic section then:
$$
\<\psi,\Psi\>\,
=\,n\big(2\b(\a^2-\b^2)\,-\,\xi.(\b^2)\big)\eta
$$
Hence if $\b=0$ then every harmonic section is a
harmonic map.  If $\a=0$ then a harmonic section is a
harmonic map if and only if:
$$
\b(\b^2+(\xi.\b))\,=\,0
$$ 
Thus either $\b=0$ or $\xi.\b=-\b^2$ (pointwise).  Therefore $\s$  is a
harmonic map if and only if $\nabla\b=-\b^2\xi$.  If $M$ is compact then
$\b$ achieves a maximum (resp\. minimum) at $x_1$ (resp\. $x_2$), say. 
Therefore $\nabla\b(x_1)=0=\nabla\b(x_2)$, hence $\b(x_1)=0=\b(x_2)$, so
$\b=0$.  The characterization of harmonic maps when $\dim
M\geqs5$ follows from Theorem 4.7.
\qed
\enddemo
\remark{Remarks}
\flushpar
(1)\quad
For a Sasakian structure, $\s$ is a harmonic map.  This generalizes Example
4.5, and strengthens the result of \cite{25} that the characteristic field
on a Sasakian manifold is harmonic.  For a Kenmotsu manifold, $\s$ is a
harmonic section, but never a harmonic map.  For a compact $\b$-Kenmotsu
manifold,
$\s$ is never a harmonic map.
\flushpar
(2)\quad
Another way of phrasing Theorem 4.8 is that harmonic sections are
characterized by $\Ric^*(\xi)=0$, and harmonic maps are characterized by
$\,\Ric^*(\xi)=0=\Ric(\xi)$, where the {\sl \midstar Ricci curvature\/} of
an almost contact manifold is defined:
$$
\Ric^*(X,Y)\,=\,\R(\Phi)(X,\phi Y).
$$
\endremark
\bigskip
\head 
5. Hypersurfaces
\endhead 
We continue the discussion from \S4 of an oriented hypersurface $M$
of an almost Hermitian manifold $(\tilde M,g,J)$, with the induced almost
contact structure (4-1).  Let $A\colon TM\to TM$ (resp\. $\a$) be the shape
operator (resp\. second fundamental form):
$$
AX\,=\,-\tilde\nab X\nu,
\qquad
\a(X,Y)\,=\,\<X,AY\>
\tag5-1
$$ 
\proclaim{5.1 Lemma}
For all $X\in TM$ we have:
\smallskip\noindent
\rom{(a)}\quad
$\tilde\nab XJ(\xi)\,
=\,\phi^2AX\,-\,J\,\nab X\xi$\rom;
\quad\qquad
\rom{(b)}\quad
$\tilde\nab XJ(\nu)\,
=\,\phi AX\,-\,\nab X\xi$
\endproclaim
\demo{Proof}
First note that $\tilde\nab XJ(\nu)$ is orthogonal to the holomorphic
$2$-plane $\xi\wedge\nu$. 
For any $Z$ in $\D$ it follows from (5-1) and (4-1) that:
$$
\<\tilde\nab XJ(\nu),Z\>\,
=\,\<-\tilde\nab X\xi-J\,\tilde\nab X\nu,Z\>\,
=\,\<\phi AX-\nab X\xi,Z\>
$$
This establishes (b), from which (a) follows since $\xi=-J\nu$.
\qed
\enddemo
We introduce the following skew-symmetric $(1,1)$ tensor $\Gamma_1$ on
$M$:
$$
\Gamma_1(X)\,=\,\tilde\nab\nu J(X)\,
-\,\<\tilde\nab\nu J(X),\nu\>\nu,
\quad\text{for all $X\in TM$.}
\tag5-2
$$
\proclaim{5.2 Proposition}
For all $X$ in $TM$ and $Z$ in $\D$ the components of
$\tilde\nabsq XXJ(Z)$ are:
\TagsOnLeft
$$
\align
\bar\nabsq XXJ(Z)\,
&+\,\a(X,X)\,\phi^2\'\Gamma_1 Z\,
+\,r(\nab X\xi,J\,\nab X\xi)Z
\tag"\rom{(a)}\quad In $\D$:" \\
&-\,r(\phi^2AX,\phi AX)Z\,+\,2\'r(\phi^2AX,\nab Y\xi)Z
\endalign
$$
\TagsOnRight
\rom{(b)}\quad
In the $\xi$-direction:
$$
\align
\nab X\a(X,Z)\,
&-\,\<\nabsq XX\xi,JZ\>\,
+\,\a(X,X)\<\Gamma_1\xi,Z\>\,
+\,2\<(\bar\nabla J\otimes\nabla\xi)(X,X),\,Z\> \\
&-\,\a(X,\xi)\a(X,JZ)\,
-\,2\'\a(X,\xi)\<\nab X\xi,Z\> 
\endalign
$$
\rom{(c)}\quad
In the $\nu$-direction:
$$
\align
\nab X\a(X,JZ)\,
&+\,\<\nabsq XX\xi,Z\>\,
+\,\a(X,X)\<\tilde\nab\nu J(\nu),Z\>\,
+\,2\'\a(\bar\nab XJ(Z),X) \\
&+\,\a(X,\xi)\a(Y,Z)\,
-\,2\'\a(X,\xi)\<\nab X\xi,JZ\>
\endalign
$$
\endproclaim
\demo{Proof}
If $Z$ is extended to a vector field then:
$$
\tilde\nabsq XXJ(Z)\,
=\,[\'\tilde\nabsq XX,J\']Z\,
-\,2\,\tilde\nab XJ(\tilde\nab XZ)
\tag5-3
$$ 
For convenience, extend $X\in T_xM$ to a vector field on $M$ with $\nabla
X(x)=0$, and $Z$ to a section of $\D$ with $\bar\nabla Z(x)=0$.  Then:
$$
\tilde\nab XX\,=\,\a(X,X)\nu,\qquad
\tilde\nab XZ\,=\,\a(X,Z)\nu\,-\,\<\nab X\xi,Z\>\xi
$$
It follows that:
$$
\align
[\'\tilde\nabsq XX,J\']Z\,
&=\,\tilde\nab X\tilde\nab X(JZ)\,-\,\tilde\nab{\tilde\nab XX}(JZ)\,
-\,J\big(\tilde\nab X\tilde\nab XZ\,-\,\tilde\nab{\tilde\nab XX}Z\big) \\
&=\,[\'\tilde\nab X\tilde\nab X,J\']Z\,
-\a(X,X)\,\tilde\nab\nu J(Z)
\endalign
$$
Routine calculations utilizing (4-1) and (5-1) then yield:
$$
\gather
\aligned
\tilde\nab X\tilde\nab X(JZ)\,
&=\,\tilde\nab X\big(\bar\nab X(JZ)\,
-\,\<JZ,\nab X\xi\>\xi\,
+\,\a(X,JZ)\nu\big) \\
&=\,\bar\nabsq XX(JZ)\,
-\,\<\nab X\xi,JZ\>\nab X\xi\,
-\,\a(X,JZ)AX \\
&\quad
-\,\<\nabsq XX\xi,JZ\>\xi\,
+\,2\<\bar\nab XJ(\nab X\xi),Z\>\xi \\
&\qquad
+\,\nab X\a(X,JZ)\nu\,
+\,2\'\a(X,\bar\nab XJ(Z))\nu\,
-\,2\'\a(X,\xi)\<\nab X\xi,JZ\>\nu, 
\endaligned
\vspace{1ex}
\aligned
J\,\tilde\nab X\tilde\nab XZ\,
&=\,J\,\tilde\nab X\big(\bar\nab XZ\,
-\,\<Z,\nab X\xi\>\xi\,
+\,\a(X,Z)\nu\big) \\
&=\,J\,\bar\nabsq XXZ\,
-\,\<\nab X\xi,Z\>J\,\nab X\xi\,
-\,\a(X,Z)\phi AX \\
&\quad
-\,\nab X\a(X,Z)\xi\,
+\,2\'\a(X,\xi)\<\nab X\xi,Z\>\xi\,
-\,\<\nabsq XX\xi,Z\>\nu\,
-\,\a(X,\xi)\a(X,Z)\nu
\endaligned
\endgather
$$
The terms in (5-3) involving $\tilde\nabla J$ may be evaluated using
Lemma 5.1, and the components of $\tilde\nabla^2J$ extracted.  For
the $\D$-component, it follows from (5-3) with $\bar\nabla$ in place of
$\tilde\nabla$ that:
$$
\bar\nabsq XX(JZ)\,-\,J\,\bar\nabsq XXZ\,
=\,\bar\nabsq XXJ(Z)
$$
The $\D$-component then simplifies using (3-5).
\qed
\enddemo
We denote by $H$ the (scalar) mean curvature of $M$:
$$
(2n+1)H\,=\,\a(E_i,E_i),
$$
and introduce a second skew-symmetric $(1,1)$ tensor $\Gamma_2$ on $M$:
$$
\Gamma_2(X)\,=\,\tilde\nabsq\nu\nu J(X)\,
-\,\<\tilde\nabsq\nu\nu J(X),\nu\>\nu,
\quad\text{for all $X\in TM$.}
\tag5-4
$$
We also abbreviate:
$$
\<\bar\nabla J,A\>(Z)\,=\,\<\bar\nab{E_i}J(Z),AE_i\>,
\qquad
\tilde\tau(J)\,=\,\tfrac14\,[\'\tilde\nabla^*\tilde\nabla J,J\']
$$
\proclaim{5.3 Proposition}
For all $Z,W$ in $\D$ the non-zero components of $\tilde\tau(J)|_M$ are
given by:
\TagsOnLeft
$$
\align
2\<\tilde\tau(J)Z,W\>\,
&=\,2\<\bar\tau(J)Z,W\>\,
-\,\tfrac12\<[\'\Gamma_2,\phi\']Z,W\>\,
+\,(2n+1)H\<\Gamma_1Z,JW\>
\tag a \\
&\qquad
+\,\<J\,\nab{AZ}\xi-\nab{AJZ}\xi,W\>\,
-\,\<J\,\nab{AW}\xi-\nab{AJW}\xi,Z\> \\
\vspace{2ex}
2\<\tilde\tau(J)Z,\xi\>\,
&=\,\nabla^*\a(JZ)\,
-\,\<\bar\nabla J,A\>(Z)\,
+\,\<\nabla^*\nabla\xi+J\'T(\phi),Z\>
\tag b \\
&\qquad
-\,\<A\xi,AZ\>\,
+\,2\<\nab{A\xi}\xi,JZ\>
\endalign
$$
\endproclaim
\TagsOnRight
\demo{Proof}
Since $r(E,JE)$ commutes with $J$, Proposition 5.2\,(a) implies:
$$
\align
4\<\tilde\tau(J)Z,W\>\,
&=\,4\<\bar\tau(J)Z,W\>\,
-\,\<[\'\tilde\nabsq\nu\nu J,J\']Z,W\>\,
-\,(2n+1)H\<[\'\phi^2\Gamma_1,J\']Z,W\> \\
&\qquad
+\,2\sum_i\<[\'r(AE_i,\nab{E_i}\xi),J\']Z,W\>, \\
\endalign
$$
which expands as stated, noting that:
$$
\gather
\<[\'\phi^2\Gamma_1,J\']Z,W\>\,
=\,-\<[\'\Gamma_1,\phi\']Z,W\>\,
=\,-2\<\Gamma_1Z,JW\>, \\
\<[\'\tilde\nabsq\nu\nu J,J\']Z,W\>\,
=\,\<\Gamma_2(JZ),W\>\,+\,\<\Gamma_2Z,JW\>\,
=\,\<[\'\Gamma_2,\phi\']Z,W\>
\endgather
$$
On the other hand since $J\xi=\nu$ we have:
$$
4\<\tilde\tau(J)Z,\xi\>\,
=\,\<\tilde\nabla^*\tilde\nabla J(JZ),\xi\>\,
+\,\<\tilde\nabla^*\tilde\nabla J(Z),\nu\>,
$$
and by Proposition 5.2\,(b) and (c):
$$
\align
\<\tilde\nabla^*\tilde\nabla J(JZ),\xi\>\,
+\,2\<T(\phi),JZ\>\,
&=\,\nabla^*\a(JZ)\,
+\,\<\nabla^*\nabla\xi,Z\>\,
-\,\<A\xi,AZ\>\,
+\,2\<\nab{A\xi}\xi,JZ\> \\
&=\,\<\tilde\nabla^*\tilde\nabla J(Z),\nu\>\,
+\,2\<\bar\nabla J,A\>(Z)
\endalign
$$
Finally note that (a) and (b) describe all the non-zero components of
$\tilde\tau(J)$, because:
$$
\<\tilde\tau(J)Z,\nu\>\,=\,\<\tilde\tau(J)(JZ),\xi\>
\quad\text{and}\quad
\<\tilde\tau(J)\xi,\nu\>\,
=\,0.
\tag"\qed"
$$
\enddemo
Our first application of Proposition 5.3 is when $\tilde M=M\times\Bbb R$. 
Then $\tilde M$ acquires a canonical almost Hermitian structure $J$,
by applying (4-1) to each leaf of the canonical hypersurface foliation, with
$\nu=d/dt$.  Let $\tilde\s$ be the section of the twistor bundle over
$\tilde M$ parametrizing $J$.
\proclaim{5.4 Theorem}
The almost Hermitian structure on $M\times\Bbb R$ is harmonic if
and only if:
$$
\bar\tau(J)\,=\,0
\qquad\text{and}\qquad
\nabla^*\nabla\xi\,=\,|\'\nabla\xi\'|^2\xi\,-\,J\'T(\phi)
$$ 
Furthermore:
$$
\<\tilde\psi,\tilde\Psi\>\,=\,\<\psi_1,\Psi_1\>\,
+\,\tfrac12\<\psi_2,\Psi_2\>
$$
If any two of $\xi$, $\s$, $\tilde\s$ are harmonic sections
(resp\. maps) then so is the third.
\endproclaim
\demo{Proof}
We have $\a=0$, and $\Gamma_1=0=\Gamma_2$, by the geometry of the
product metric.  Proposition 5.3 therefore reduces to:
$$
\align
\<\tilde\tau(J)Z,W\>\,
&=\,\<\bar\tau(J)Z,W\>,
\qquad
2\<\tilde\tau(J)Z,\xi\>\,
=\,\<\nabla^*\nabla\xi+J\'T(\phi),Z\>,
\endalign
$$
and since the $\xi$-component of $\nabla^*\nabla\xi$ is $|\'\nabla\xi\'|^2$
it follows that $\tilde\tau(J)=0$ if and only if:
$$
\bar\tau(J)=0
\qquad\text{and}\qquad
\nabla^*\nabla\xi\,-\,|\'\nabla\xi\'|^2\xi\,+\,J\'T(\phi)\,=\,0
$$
\par
Let $\{F_i\}$ be a local orthonormal frame of $\D$, extended to a local
frame $\{E_i\}$ of $TM$ by defining $E_{2n+1}=\xi$, and a local
frame $\{\tilde E_i\}$ of $T\tilde M$ by defining $\tilde E_{2n+2}=\nu$.
Since $\tilde R$ vanishes when applied to $\nu$ (in any of its arguments),
for all $X\in TM$ we have:
$$
\align
4\<\tilde\psi,\tilde\Psi\>(X)\,
&=\,\<\tilde R(\tilde E_i,X)\tilde E_j,\,
J\,\tilde\nab{\tilde E_i}J(\tilde E_j)\>\,
=\,\<R(E_i,X)E_j,\,J\,\tilde\nab{E_i}J(E_j)\> \\
&=\,\<R(E_i,X)\xi,\nab{E_i}\xi\>\,
+\,\<R(E_i,X)F_j,\,J\,\tilde\nab{E_i}J(F_j)\>,
\quad\text{by Lemma 5.2,} \\
&=\,\<\psi_2,\Psi_2\>(X)\,
+\,4\<\psi_1,\Psi_1\>(X) \\
&\qquad
-\,\big\<R(E_i,X)F_j,\,\<\tilde\nab{E_i}J(F_j),\nu\>\xi\big\>,
\quad\text{by Theorem 3.4 and Lemma 4.2,} \\
&=\,4\<\psi_1,\Psi_1\>(X)\,+\,2\<\psi_2,\Psi_2\>(X),
\quad\text{by Lemma 5.2.}
\endalign
$$
The relationships between the harmonicity of $\xi$, $\s$ and $\tilde\s$
follow by comparison with Theorems 3.2 and 3.4.
\qed
\enddemo
The following consequence of Theorems 5.4 and 4.8 generalizes the
result of \cite{27} that the almost Hermitian structure on the Hopf manifold
$S^{2n+1}\times S^1$ is harmonic.
\proclaim{5.5 Corollary}
Suppose $M$ is a Sasakian manifold.  Then the section
$\tilde\s$ parametrizing the canonical almost Hermitian structure on
$M\times\Bbb R$ \rom(or $M\times S^1$\rom) is a harmonic map.
\endproclaim
Recall that an almost contact manifold is said to be {\sl nearly
cosymplectic\/} \cite{2} if: 
$$
\nab X\phi(X)=0,
\qquad\text{for all $X\in TM$.}
\tag5-5
$$
\proclaim{5.6 Corollary}
Suppose $M$ is a nearly cosymplectic manifold with parallel
characteristic field.  Then $\s$ is a harmonic map. 
\endproclaim
\demo{Proof}
It follows from (5-5) that $\tilde\nab XJ(X)=0$ for all $X\in TM$. 
Moreover:
$$
\tilde\nab\nu J(X)\,=\Gamma_1(X)\,=\,0,\qquad
\tilde\nab XJ(\nu)\,=\,-\nab X\xi,
$$
so $\tilde M$ is nearly K\"ahler precisely when $\xi$ is parallel.  Now
$\tilde\s$ is a harmonic map \cite{29}, and therefore $\s$ is a harmonic
map by Theorem 5.4.  
\qed
\enddemo
By analogy with \cite{14} we say a nearly cosymplectic structure is
{\sl strict\/} if $\nab X\phi\neq0$ for all non-zero $X$; for example,
if $\phi$ is induced by a strict nearly K\"ahler structure on $\tilde M$. 
Corollary 5.6 does not apply if $\phi$ is strict, because nearly
cosymplectic manifolds satisfy:
$$
\nab\xi\phi(X)\,=\,-\phi\,\nab X\xi
$$
We therefore consider the
case when $\tilde M$ is nearly K\"ahler.  We may (locally) extend $\nu$ to
a vector field on a neighbourhood of $M$ in $\tilde M$ by parallel
translation along the normal geodesics.  Any vector field on $M$ may be
extended in a similar way.  Since $\tilde M$ is nearly K\"ahler its
geodesics are holomorphically planar.  In particular, the holomorphic
$2$-plane $\xi\wedge\nu$ remains holomorphic when parallel transported along
a $\nu$-geodesic, and hence:
$$
\tilde\nab\nu J(\nu)\,
=\,-\tilde\nab\nu\xi\,-\,J\,\tilde\nab\nu\nu\,
=\,0
\tag5-6
$$
Consequently, by Lemma 5.1\,(b), the definition (5-2) of $\Gamma_1$
simplifies to:
$$
\Gamma_1X\,=\,\tilde\nab\nu J(X)\,
=\,\nab X\xi\,-\,\phi AX
\tag5-7
$$
Furthermore $\Gamma_1$ is a section of $\fM_1$ (see
\S3), for it follows from (4-1) that:
$$
\{\Gamma_1,\phi\}X\,
=\,-\<X,\tilde\nab\nu J(\nu)\>\xi\,
-\,\eta(X)\,\tilde\nab\nu J(\nu)\,=\,0.
$$
\proclaim{5.7 Proposition}
Suppose $\tilde M$ is nearly K\"ahler.
\flushpar
\rom{(a)}\quad 
The integral curves of $\xi$ are geodesics if and only if $\xi$ is a
principal direction.
\flushpar
\rom{(b)}\quad
The characteristic field $\xi$ is Killing if and only if $\;[\'A,\phi\']=0$.
\flushpar
\rom{(c)}\quad
$M$ inherits a contact metric structure if and only if $\;\Gamma_1=0$ and 
$\;\{A,\phi\}=-2\phi$.
\flushpar
\rom{(d)}\quad
$M$ inherits a nearly cosymplectic structure if and only if 
$\;\a=(2n+1)H\eta\otimes\eta$.
\endproclaim
\par\pagebreak
\demo{Proof}
\flushpar
(a)\quad
Since $\Gamma_1\xi=0$ it follows from (5-7) that $\;\nab\xi\xi=\phi
A\xi$.
\flushpar
(b)\quad
It follows from (5-7) and the skew-symmetry of $\Gamma_1$ that:
$$
\<\nab X\xi,Y\>\,+\,\<X,\nab Y\xi\>\,
=\,\<[\'\phi,A\']X,Y\>
$$
\flushpar
(c)\quad
From (5-7):
$$
\align
2\'d\eta(X,Y)\,
&=\,\nab X\eta(Y)\,-\,\nab Y\eta(X)\,
=\,\<\nab X\xi,Y\>\,-\,\<X,\nab Y\xi\> \\
&=\,-\<X,(2\Gamma_1+\{A,\phi\})Y\>
\tag5-8
\endalign
$$
Therefore if $\Gamma_1=0$ and $\{A,\phi\}=-2\phi$ then $\phi$ is a
contact metric structure.  Conversely, since contact metric structures
satisfy $\nab\xi\phi=0$ it follows from (5-7), Lemma 4.2 and (4-2) that
for all $Z\in\D$:
$$
\Gamma_1Z\,
=\,\tilde\nab\nu J(Z)\,
=\,-\tilde\nab\xi J(JZ)\,
=\,-\bar\nab\xi J(JZ)\,
=\,-\nab\xi\phi(JZ)\,=\,0
\tag5-9
$$
Hence $\Gamma_1=0$.  It then follows from (5-8) that $\{A,\phi\}=-2\phi$.
\flushpar
(d)\quad
Proved in \cite{3, Theorem 6.13}.
\qed
\enddemo
Proposition 5.7\,(c) generalizes the criterion \cite{20} for contact
metric hypersurfaces of K\"ahler manifolds.  In particular, if $\tilde M$
is a strict nearly K\"ahler manifold then no contact metric hypersurfaces
exist.
\par
Concerning $\Gamma_2$, it follows from (5-6) that:
$$
\tilde\nabsq\nu\nu J(\nu)\,
=\,\tilde\nab\nu\tilde\nab\nu J(\nu)\,
=\,\tilde\nab\nu(\tilde\nab\nu J(\nu))\,
-\,\tilde\nab\nu J(\tilde\nab\nu\nu)\,
=\,0,
$$
so (5-4) simplifies to:
$$
\Gamma_2(X)\,=\,\tilde\nabsq\nu\nu J(X),
\quad\text{for all $X\in TM.$}
$$
Since $\tilde M$ is nearly K\"ahler, it follows from 
\cite{13, Proposition 2.3} that:
$$
\<\tilde\nabsq\nu\nu J(-),J(-)\>\,
=\,-\<\tilde\nab\nu J(-),\tilde\nab\nu J(-)\>,
$$
and hence, for all $X,Y\in TM$:
$$
\<\Gamma_2X,Y\>\,=\,\<\phi\Gamma_1X,\Gamma_1Y\>
$$
Therefore, since $\Gamma_1$ is a section of $\fM_1$,
$\Gamma_2$ is a section of $\fH$. 
\proclaim{5.8 Proposition}
If $\tilde M$ is nearly K\"ahler then for all $Z,W\in\D$ we have:
\smallskip\noindent
\rom{(a)}\qquad
$\<\bar\tau(J)Z,W\>\,=\,2(2n+1)H\<J\,\Gamma_1Z,W\>\,
-\,2\<[\{A,\phi\},\Gamma_1\']Z,W\>$
\smallskip\noindent
\rom{(b)}\qquad
$\<\nabla^*\nabla\xi+JT(\phi),Z\>\,
=\,\<(\Gamma_1J-A)Z,A\xi\>\,
-\,\nabla^*\a(JZ)$
\endproclaim
\demo{Proof}
Since $\tilde M$ is nearly K\"ahler we have $\tilde\tau(J)=0$
\cite{27}.
\flushpar
(a)\quad 
Since $\Gamma_2$ commutes with $\phi$, Proposition 5.3\,(a) reduces to:
$$
-\tfrac12\<\bar\tau(J)Z,W\>\,
=\,(2n+1)H\<\Gamma_1Z,JW\>\,+\,\Sigma,
$$
where:
$$
\Sigma\,=\,\<\nab{AW}\xi,JZ\>\,
-\,\<\nab{AZ}\xi,JW\>\,
+\,\<\nab{AJW}\xi,Z\>\,-\,\<\nab{AJZ}\xi,W\>
$$
By (5-7), noting the cancellation of terms involving $A^2$, we obtain:
$$
\align
\Sigma\,
&=\,\<\Gamma_1AW,JZ\>\,
-\,\<\Gamma_1AZ,JW\>\,
+\,\<\Gamma_1AJW,Z\>\,
-\,\<\Gamma_1AJZ,W\> \\
&=\,\<(\phi\'\Gamma_1\'A\,-\,\Gamma_1\'A\'\phi\,
+\,\phi A\'\Gamma_1\,-\,A\'\Gamma_1\'\phi)Z,W\>,
\quad\text{by (4-1)} \\
&=\,\<[\'\{A,\phi\},\Gamma_1\']Z,W\>,
\quad\text{since $\Gamma_1$ anticommutes with $\phi$.}
\endalign
$$
(b)\quad
Define a symmetric tensor $S\colon\D\to\D$ by $S=\phi^2A\'|\'\D$.  Since
$\tilde M$ is nearly K\"ahler:
$$
\<\bar\nabla J(Z),A\>\,
=\,\<\bar\nab ZJ,S\>\,
+\,\<\bar\nab\xi J(Z),A\xi\>\,
=\,\<J\Gamma_1Z,A\xi\>,
$$
by the skew-symmetry of $\bar\nab ZJ\colon\D\to\D$, and (5-9).
Furthermore, by (5-7):
$$
\<\nab{A\xi}\xi,JZ\>\,
=\,\<(A-\Gamma_1J)Z,A\xi\>,
$$
and Proposition 5.3\,(b) simplifies as stated.
\qed
\enddemo
In view of Proposition 5.7 it is reasonable to consider hypersurfaces with
$\Gamma_1=0$ and $\xi$ a principal direction, which include all contact
metric hypersurfaces; in fact, scrutiny of the proof shows they are
characterized by the intrinsic condition $\nab\xi\phi=0$.  We say that 
{\sl $\phi$ is harmonic\/} if $\s$ is a harmonic section.
\proclaim{5.9 Theorem}
Let $M$ be a hypersurface of a nearly K\"ahler manifold $\tilde M$, with
$\Gamma_1=0$ and $\xi$ a principal direction.  The induced almost contact
structure $\phi$ is harmonic if and only if $\xi$ is harmonic, if and only
if $\nabla^*\a$ is characteristic.  If $\nu$
\rom(equivalently, $\xi$\rom) is a Ricci-principal direction in $\tilde M$
then $\phi$ is harmonic precisely when $\nabla H$ is characteristic.
Thus if $\tilde M$ is Einstein and $M$ has constant mean curvature then
$\phi$ is harmonic.
\endproclaim
\demo{Proof}
Since $\xi$ is a principal direction, we may pick an orthonormal frame
$\{F_i\}$ of $\D$ with $AF_i=\kappa_iF_i$, say.  Then since $\bar\nab\xi
J=0$ by (5-9), it follows from (5-7) that:
$$
T(\phi)\,=\,\bar\nab{F_i}J(\nab{F_i}\xi)\,
=\,\bar\nab{F_i}J(\phi AF_i)\,
=\,-\kappa_i\,J\,\bar\nab{F_i}J(F_i)\,
=\,0,
$$
by Lemma 4.2.  Furthermore
$\bar\tau(J)=0$ by Proposition 5.8\,(a).  It therefore follows from
Theorem 3.2 that $\phi$ is harmonic if and only if $\xi$ is harmonic. 
But by Proposition 5.8\,(b) this is the case precisely when
$\nabla^*\a(\D)=0$.
\par
It follows from Codazzi's equation:
$$
\nab X\a(Y,Z)\,-\,\nab Y\a(X,Z)\,=\,\<\tilde R(X,Y)Z,\nu\>,
$$
that:
$$
\nabla^*\a(Z)\,
=\,\widetilde{\Ric}(Z,\nu)\,-\,(2n+1)dH(Z)\,
=\,-\widetilde{\Ric}(Z,\xi)\,-\,(2n+1)dH(Z),
$$
since the Ricci curvature of a nearly K\"ahler manifold is $J$-invariant
\cite{18}.  The $J$-invariance also implies $\widetilde{\Ric}(\nu,\xi)=0$,
and therefore $\nu$ (equivalently, $\xi$) is an eigenvector of
$\widetilde{\Ric}$ if and only if $\widetilde{\Ric}(\nu,\D)=0$
(equivalently, $\widetilde{\Ric}(\xi,\D)=0$).  Therefore $\nabla^*\a$ is
characteristic if and only if $dH$ is characteristic.
\qed
\enddemo
\proclaim{5.10 Theorem}
Suppose $M$ is a totally umbilical hypersurface of a nearly K\"ahler
manifold $\tilde M$.   If $M$ is totally geodesic and $\xi$ is harmonic
then the induced almost contact metric structure is
harmonic.  The converse is true if $\tilde M$ is strict.
\endproclaim
\demo{Proof}
Since $\a=Hg$, Proposition 5.8 reduces to:
$$
\gather
\bar\tau(J)\,=\,2(2n-3)HJ\'\Gamma_1 
\tag5-10 \\
\nabla^*\nabla\xi\,
-\,|\'\nabla\xi\'|^2\xi\,
+\,J\'T(\phi)\,
=\,-\phi\,\nabla H
\tag5-11
\endgather
$$
Therefore if $H=0$ and $\xi$ is harmonic then $\bar\tau(J)=0$ and
$T(\phi)=0$, so $\phi$ is harmonic by Theorem 3.2.  Conversely if
$\tilde M$ is strict and $\phi$ is harmonic then (5-10) implies $H=0$, and
comparison of (5-11) with Theorem 3.2 implies $\xi$ is harmonic.
\qed
\enddemo
It follows from Theorem 5.10 that the only (round) hypersphere of the
nearly K\"ahler $6$-sphere whose induced almost contact structure is
harmonic is the equator, which is nearly cosymplectic.  In
particular, the almost contact structure induced on the hypersphere of
radius $1/\sqrt2$, which is {\sl nearly Sasakian\/} \cite{5}, is not
harmonic. Interestingly, the characteristic vector field of any round
hypersphere is Killing (eg\. Proposition 5.7), and hence harmonic (cf\.
Theorem 4.3).  Regarding the nearly cosymplectic $5$-sphere:
\proclaim{Theorem 5.11}
The section parametrizing the standard nearly cosymplectic structure on
$S^5$ is a harmonic map.
\endproclaim
\demo{Proof}
Let $\{F_i\}$ be an orthornormal frame in $\D$, extended to a frame
$\{E_i\}$ of $TM$ with $E_{2n+1}=\xi$.  
Note that (5-5) implies $\bar\nabla J(F,F)=0$ for all $F\in\D$.
From Theorem 3.4:
$$
\align
4\<\psi_1,\Psi_1\>(X)\,
&=\,\<R(E_i,X)F_j,\,J\,\bar\nab{E_i}J(F_j)\>\,
=\,\<\<X,F_j\>E_i-\<E_i,F_j\>X,\,J\,\bar\nab{E_i}J(F_j)\> \\
&=\,-2\<J\,\bar\nab{F_i}J(F_i),X\>\,=\,0, \\
\<\psi_2,\Psi_2\>(X)&=\,\<R(E_i,X)\xi,\nab{E_i}\xi\>\,
=\,\<X,\xi\>\div\xi\,-\,\<X,\nab\xi\xi\>\,=\,0.
\tag"\qed"
\endalign
$$
\enddemo

\enddocument
\end